\documentclass[letterpaper, 10 pt, peerreview, onecolumn]{ieeeconf}  

\IEEEoverridecommandlockouts                              
\usepackage{graphics} 
\usepackage{epsfig} 
\usepackage{mathptmx} 
\usepackage{times} 
\usepackage{amsmath} 
\usepackage{amssymb}  
\usepackage{esvect}

\newtheorem{theorem}{Theorem}
\newtheorem{remark}{Remark}
\newtheorem{assumption}{Assumption}
\newtheorem{definition}{Definition}
\newtheorem{proposition}{Proposition}

\title{\LARGE \bf
A barrier function approach to constrained Pontryagin-based Nonlinear Model Predictive Control
}

\author{Michele Pagone, Mattia Boggio, Carlo Novara, Anton Proskurnikov, Giuseppe C. Calafiore
\thanks{The authors are with the Department of Electronics and Telecommunications,
        Corso Duca degli Abruzzi, 24, 10129 Torino, Italy.
        {\tt\small \{michele.pagone, mattia.boggio, carlo.novara, anton.proskurnikov, giuseppe.calafiore\}@polito.it}}%
}

\begin{document}

\maketitle
\thispagestyle{empty}
\pagestyle{empty}

\begin{abstract}

A Pontryagin-based approach to solve a class of constrained Nonlinear Model Predictive Control problems is proposed which employs the method of barrier functions for dealing with the state constraints. Unlike the existing works in literature the proposed method is able to cope with nonlinear input and state constraints without any significant modification of the optimization algorithm. A stability analysis of the closed-loop system is carried out by using the $\mathcal{L}_2$-norm of the predicted state tracking error as a Lyapunov function. Theoretical results are tested and confirmed by numerical simulations on the Lotka-Volterra prey/predator system.

\end{abstract}

\section{INTRODUCTION}\label{sec1}
Over the last years, Model Predictive Control (MPC) has been accepted as a powerful control tool for a wide range of technological  applications~\cite{richalet1978,qin2000}, thanks to its capability to design control algorithms for multivariate systems under state, input, and output constraints. The resulting controller also provide optimality of a predefined performance index. 

The key point of the MPC design is the method for addressing optimal control problems (OCP) with receding horizon. To cope with nonlinear dynamics and constraints, as well as with non-convex performance indexes, Nonlinear MPC (NMPC) have been introduced (see, e.g.~\cite{diehl2007} and references therein). To find the global optimum in this situation is difficult, optimization algorithms are computationally intensive and, in general, the solution rarely admits an explicit closed-form representation~\cite{allgower2004, cisneros2020}.

In this paper, we propose a solution that is based on the Pontryagin's Minimum Principle (PMP)\cite{pontryagin1962}: under some assumptions on the Hamiltonian function, we can obtain an explicit control law - as function of the state and the co-state - even if the system dynamics and/or constraints are nonlinear. The price paid for this is the necessity to solve a Two-Points Boundary Value Problem (TPBVP) in order to find the state and co-state functions. The first applications of the PMP to receding horizon control date back to works by \cite{ohtsuka1993, ohtsuka1997} and \cite{palanki1993} who have also established important higher degree optimality conditions based on the theory of Lie algebras.

Although TPBVP problems usually cannot be solved analytically, a number of efficient numerical algorithms to solve OCP in real time have been proposed \cite{srinivasan2003} such as, e.g., the stabilized continuation method \cite{ohtsuka1993} and its accelerated versions \cite{ohtsuka2003}, the Newton-type algorithm \cite{deng2019} and the extended modal series method, approximating OCP with nonlinear constraints by standard LQR problems \cite{devia2018}. An efficient \emph{active set} method of solving discrete-time PMP equations arising in MPC problems with input and terminal state constraints was developed in \cite{cannon2008}. Continuous-time OCP can be accurately approximated by discrete-time ones as demonstrated by the recent work \cite{dontchev2020}.

Whereas initial and terminal state constraints can be accommodated by existing PMP-based MPC algorithms, direct application of PMP becomes problematic in the situation where the state vector is constrained at any time 
\cite{diehl2007,malisani2016}. In this situation, the differential equations of PMP are different for constrained and unconstrained trajectories. It arises the necessity of `tailoring' unconstrained and constrained pieces of the trajectory by imposing additional interior tangency conditions at the junctions points \cite{pesch1994, bonnard2003}. This substantially complicates the solution of TPBVP in real time except for the situations where the optimal solution structure is known a priori.

An alternative way to cope with state or mixed input-state constraints is based on the use of barrier functions that arise as penalty terms in the objective function. A general methodology to get rid of relaxing both state and input constraints by introducing penalty terms has been proposed in \cite{graichen2009} under the assumption that the nonlinear system has a well-defined relative degree. A similar approach has been proposed for a special type of constraints in \cite{kovaltchouk2015}. In this paper, we further elaborate the approach proposed in the example from \cite{ohtsuka1997}, where the state constraint is replaced by an appropriate penalty term in the cost functional, without significantly modify the algorithm of solving OCP compared to the unconstrained case. Unlike \cite{graichen2009,kovaltchouk2015}, input constraints does not need to be relaxed and can be tackled by the standard PMP. 

The penalty function method proposed in this paper is concerned with defining a methodology for accounting the state constraints within the TPBVP, without affecting the differential equation solution feasibility. This latter aspect was widely discussed by \cite{hauser2006}, which pointed out that, when employing the classical log-barrier function, some TPBVP feasibility issues can arise. A similar approach can be found \cite{suwartadi2010}: a Lagrangian-barrier function based method which adds the state constraints as a logarithmic term to the objective function. As remarked in \cite{malisani2016}, the penalty functions methods can be divided into two different classes: exterior and interior. We focus on the interior penalty methods since they generate only feasible solutions. This can be an interesting particularity in numerous nonlinear and non-convex applications: satisfaction of constraints is more important than optimality (see also \cite{murray2003} and the reference therein).

We propose a class of Gaussian-like penalty function. Thanks to this approach, the solution to the system of differential equation is (almost) always guaranteed. An important advantage of the proposed penalty methodology relies in the relaxation on the constraints and penalty function assumptions, in particular, the penalty function has to be only $C^1$-smooth unlike the approach from \cite{graichen2009}. 

To sum up, the proposed NMPC framework shows the following advantages: i) conversely to the numerical methods where a discretization of state, input and constraints before optimization is required \cite{boiroux2019}, the PMP-based solution does not need the input parametrization anymore, resulting in a better accuracy in tracking the reference; ii) the PMP-based NMPC seems to perform an more efficient trade-off between computational complexity and final reference tracking with respect to the direct methods, making him suitable for on-line applications. 

The paper is organized as follows. 
In Section \ref{NMPC_Frame} the NMPC scheme and its unconstrained Pontryagin-based solution are illustrated. The PMP-based solution of the constrained problem is shown in Section \ref{constrained}. Moreover, in Section \ref{finitegain} we propose the mathematical study about the local stability of the closed-loop system. A simulated example is presented in Section \ref{sim}. Finally, the conclusions are drawn in Section \ref{conclusion}.

\section{NMPC Framework}\label{NMPC_Frame}
Consider the following affine-in-the-input nonlinear system:
\begin{equation} \label{eq:nonlinear}
\dot{x}(t) = f(x(t)) + g(x(t))u(t)  	\\
\end{equation}
where $x \in \mathbb{R}^{n_x}$, $u \in\mathbb{R}^{n_u}$ are the state and the input, respectively. We assume that the state of system \eqref{eq:nonlinear} is measured in real time, with a sampling time $T_S$. At each time $t=t_k$, a prediction of the system state and output over the time interval $[t, t + T_p]$ is performed, where $T_p \ge T_S$ is the prediction horizon. The prediction is obtained by integrating~(\ref{eq:nonlinear}). At each time $t = t_k$, we look for an input signal $u^*(t:t+T_P)$, minimizing a suitable cost function $J\bigl(u(t:t + T_p)\bigr)$ subject to possible constraints that may occur during the system's operations. The considered NMPC cost function $J\bigl(u(t:t + T_p)\bigr)$ in the Bolza form is 
\begin{equation}\label{eq:quadratic}
\begin{split}
J &= \int_{t}^{t + T_p} \tilde{x}^T_p(\tau)\mathbf{Q}\tilde{x}_p(\tau) ~\mathrm{d}\tau~ +\\
&+  \int_{t}^{t + T_p} u^T(\tau)\mathbf{R}u(\tau)~\mathrm{d}\tau ~+ \tilde{x}^T_p(t + T_p)\mathbf{P}\tilde{x}_p(t + T_p) .
\end{split}
\end{equation}
Mathematically, at each time $t = t_k$, the following optimization problem is solved:
\begin{equation} \label{eq:functional}
\begin{split}
& u^*(t:t+T_p) = \arg \min_{u(\cdot)}J(u(t:t + T_p)) \\
& \text{subject to:} \\
& \dot{\hat{x}}(\tau) = f\bigl(\hat{x}(\tau)) + g(\hat{x}(\tau)) u(\tau)\bigr),~~\hat{x}(t) = x(t) \\
& \hat{x}(\tau) \in X_C, ~u(\tau) \in U_C, ~~\forall \tau \in [t:t+T_p] ,\\
& u(\tau)\in \mathcal{KC}\bigl([t, t + T_p]\bigr).
\end{split}
\end{equation}
$X_C$ and $U_C$ are sets describing possible constraints on the state, output and input, respectively and $\mathcal{KC}([t, t+T_p])$ is the space of piece-wise continuous functions. A receding control horizon strategy is employed: at a given time $t=t_k$, the input signal $u^*(t_k:t_k + T_p)$ is computed by solving \eqref{eq:functional}. Then, only the first optimal input value $u(t) = u^*(t_k)$ is applied to the plant, keeping it constant for $t \in [t_k, t_{k+1}]$. The remainder of the solution is discarded. Then, the complete procedure is repeated at the next time steps $t=t_{k+1}, t_{k+2}, ...$

\begin{assumption}\label{lip}
Let $f  \in \mathcal{C}^1 (\mathbb{R}^{n_u} \times \mathbb{R}^{n_x} \rightarrow \mathbb{R}^{n_x})$ and $g \mathcal{C}^1 ( \mathbb{R}^{n_u} \times \mathbb{R}^{n_x} \rightarrow \mathbb{R}^{n_x})$. 
\end{assumption}

\begin{assumption}\label{u_conv}
The admissible control set $U_C \subseteq \mathbb{R}^{n_u}$ is a a ball $U_C = \{u\in \mathbb{R}^{n_u} : \| u\|_q \leq u_{max}\}$.
\end{assumption}

\begin{assumption}
The state constraint set is  $X_C = \{x \in \mathbb{R}^{n_x} : C(x) \leq 0\}$. Here $C(x(t)) \in \mathcal{C}^1 (\mathbb{R}^{n_{x}} \rightarrow \mathbb{R})$ is, generally, a non-convex function.
\end{assumption}

\begin{remark}
The optimization problem \eqref{eq:functional} is numerically hard to tackle, since $u$ is a continuous-time signal and thus the number of decision variables is infinite. The direct solution of the OCP requires a finite parametrization of the input signal $u$ (see, e.g., \cite{boiroux2019}). For example, as illustrated in Section \ref{sim} a piece-wise constant parametrization can be assumed, with changes of value at the nodes $\tau_1,\ldots,\tau_{N}\in[t,t+T_p]$ with $N$ the number of nodes. The choice of $N>1$ can lead to satisfactory performances behaviors, but at cost of computational complexity increment. One can pick $N=1$ (corresponding to a constant input for every $\tau \in [t,t+T_P]$) in order to reduce the computational complexity of the optimization algorithm. Nevertheless, this approach could not always guarantee an acceptable level of performance. This issue is mitigated when using the PMP approach presented in the manuscript which does not require any a-priori prarametrization of the control signal. This latter does not significantly effect the algorithm computational complexity. 
\end{remark}


\subsection{Unconstrained Pontryagin-based NMPC Solution}

We neglect for the moment possible constraints on the state and the input, focusing on the case where $X_C \equiv \mathbb{R}^{n_x}$ and $U_C \equiv \mathbb{R}^{n_u}$.

According to \cite{pontryagin1962}, a necessary condition for a trajectory $x(t)$ to be the extremal path and the corresponding control $u(t)$ to be the optimal input, is that the Hamiltonian scalar function $H(x(t), u(t), \lambda(t)) \in \mathcal{C}^k (\mathbb{R}^{n_x}\times\mathbb{R}^{n_u} \times \mathbb{R}^{n_x} \rightarrow \mathbb{R})$ attains its minimum value when $u = u^*$ and while satisfying the differential equations of the dynamics in \eqref{eq:nonlinear}, the time evolution of the Lagrangian multipliers $\lambda \in \mathbb{R}^{n_x})$ (or co-state variables), and a set of boundary conditions (B.C.). The Hamiltonian is defined as
\begin{equation}\label{eq:Ham_quad}
H = \tilde{x}_p^T\mathbf{Q}\tilde{x}_p +  u^T\mathbf{R}u + \lambda^T \bigl(f + gu\bigr).
\end{equation}
The necessary conditions for optimality can be derived by analyzing the first-order variation of the augmented expression of \eqref{eq:quadratic}. The rigorous mathematical formulation of the first-order variation can be found in \cite{bryson1975}. Whereby, the Pontryagin formulation of the NMPC optimal control problem:
\begin{equation} \label{eq:functional2}
\begin{split}
&(x^*, u^*, \lambda^*) = \arg \min_{u(\cdot)}H\\
& \text{subject to:} \\
& \dot{x} = f + gu \\
& \dot{\lambda} = -\nabla_x H^T \\
& \psi = 0 \\ 
&\lambda^T_{t_i} = -\mu^T \\
&\lambda^T_{t_f} = 2\mathbf{P}\tilde{x}_{p}(t_f)
\end{split}
\end{equation}
From \eqref{eq:functional2}, we can note that the optimization problem is subject to both the state dynamics in \eqref{eq:nonlinear}, and the dynamic of the co-state variables $\lambda$, described by the so-called Euler-Lagrange differential equations. Both the state and co-state evolution must satisfy a set of boundary conditions to be imposed at the borders of the prediction horizon. The B.C. have to be satisfied by $\lambda$ and $x$ during the system evolution along the extremal path, whereas $\psi(x(t_k), x(t_f)) = 0$ are boundary conditions affecting the state at the boundaries of the prediction horizon. At each time $t = t_k$, the state value cannot be chosen arbitrarily: the continuity between two successive sampling steps must be ensured, so that $\psi =x_i - x(t_i)$. In \eqref{eq:functional2}, at $t = t_i$, $\lambda(t_i) = -\mu$ where $\mu$ is an adjoint vector whose entries can be arbitrarily picked at the beginning of the prediction horizon.

The Euler-Lagrange equations - describing the $\lambda$ time evolution - take the form of:
\begin{equation}\label{eq:EL_quad}
\begin{split}
\dot{\lambda} = -\bigl(\lambda^T\nabla_x \bigl(f(x) + g(x)u\bigr) - 2\mathbf{Q}\tilde{x}_p\bigr)^T.
\end{split}
\end{equation}
The optimal control law is obtained by minimizing the Hamiltonian with respect to $u$. By observing the Equation \eqref{eq:Ham_quad}, we have
\begin{equation}\label{eq:control_law_quad}
u^* = -\dfrac{1}{2}\mathbf{R}^{-1}\lambda^Tg(x)
\end{equation}
where $\mathbf{R}$ is constant, diagonal, positive, and invertible matrix.

By observing the PMP-based NMPC solution in \eqref{eq:functional2}, together with the optimal control law in \eqref{eq:control_law_quad}, it is clear how the optimal control problem in \eqref{eq:functional2} turns into a two-points boundary value problem. Indeed, the equations \eqref{eq:nonlinear} together with \eqref{eq:EL_quad} and the B.C. in \eqref{eq:functional2} represents a TPBVP to be solved over the prediction horizon $[t, t + T_p]$. The TPBVP solution provides the $\lambda$ and the $x$ of the explicit control laws \eqref{eq:control_law_quad}.

The TPBVP is formalized as follows:
\begin{equation}
\begin{split}
\dot{x} &= f(x) + g(x)u \\
\dot{\lambda} &= -\nabla_xH^T \\
x_i - x(t_i) &= 0 \\
\lambda^T(t_f) &= 2\mathbf{Q}\tilde{x}_p(t_f)
\end{split}
\end{equation}

\begin{remark}
Observing the optimal control laws \eqref{eq:control_law_quad}, the input $u^*(\tau)$ depends on $\lambda(\tau)$ and $x(\tau)$, whose values change at each sampling step of the TPBVP over the prediction horizon. For this reason, the PMP-based NMPC solutions does not require an a-priori parametrization of the input signal. This is a very interesting results since the OCP algorithm achieves high performances without increasing the computational complexity independently from the input parametrization.
\end{remark}


\section{Indirect Solution of the Constrained OCP}\label{constrained}
In general, the constrained case can be handled by means the indirect optimization problem only when the optimization is performed off-line, by augmenting the system with additional variables \cite{bryson1975, wang2017}. Nevertheless, when dealing with a on-line optimization process, this aspect can be tough, since it is necessary to iterate the solution in order to identify the control arcs where the constraints are active and imposing additional B.C. at the junction points.

\subsection{Input Constraints}
We consider that the input is bounded linearly, such that $U_C = \{u(t) \in \mathbb{R}^{n_u} : ~ u_{i_{min}} \leq u_i(t)\leq u_{i_{max}}, ~\forall t \}$. Consider the optimal control law \eqref{eq:control_law_quad}, for the nonlinear system \eqref{eq:nonlinear}, the optimal control $u^* \in U_C$ is:
\begin{equation}\label{con_sat}
u^* = \mathrm{sat}_{U_C}\biggl(-\dfrac{1}{2}\mathbf{R}^{-1}\lambda^Tg(x)\biggr)
\end{equation}
where the $\mathrm{sat}(\cdot)$ represents the saturation operator and it applies element-wise to the input vector. In formulae, the $i^{th}$ control component is:
\begin{equation}\label{eq:sat_con}
u_i^* = \begin{cases}
 u_{i_{min}},   &\mathrm{if}  ~-\dfrac{\lambda_ig_i(x)}{2r_i} \leq  u_{i_{min}}\\
 u_{i_{max}},  &\mathrm{if} ~-\dfrac{\lambda_ig_i(x)}{2r_i} \geq  u_{i_{max}} \ \\
-\dfrac{\lambda_ig_i(x)}{2r_i},  ~   &\mathrm{otherwise}
\end{cases}
\end{equation}
where $r_i$ is the $i^{th}$ entry of the $\mathbf{R}$ diagonal.

\begin{proposition}
For the nonlinear system \eqref{eq:nonlinear} with performance index \eqref{eq:quadratic}, if $u \in U_C$, the constrained optimal command is given by \eqref{eq:sat_con}.
\end{proposition}

\begin{proof}
From the optimal control equation we have $u^* = \arg \min_{u \in U_C}H$. For the problem at hand, since $\nabla_u \lambda^Tf(x) = 0$, we can neglect the terms not depending on the control in the Hamiltonian. Then, picking only the control-depending terms of the Hamiltonian and recalling that $\mathbf{R}$ is a diagonal positive matrix:
\begin{equation}
u^* = \arg \min_{u \in U_C}\biggl[\sum_{i=1}^{n_u} \mathbf{R}_{ii}u_i^2 + \sum_{i = 1}^{n_x}\lambda_ig_i(x)u_i\biggr].
\end{equation}
Since there are not coupled control terms, the optimal control equation can be solved by minimizing the Hamiltonian element-wise. This is straightforward, since, in this configuration, the Hamiltonian consists in a elliptic paraboloid whose main axes are parallel to the Cartesian axes. Consider the unconstrained case. Being the Hamiltonian convex with respect to $u$ we have that $H(u) \geq H(u^*) + \nabla_u H(u^*)^T(u - u^*)$ , i.e. all the admissible values of the input are enclosed in one of the halfspaces $\mathcal{H}_{++}$ delimited by the hyperplane tangent at $H$ in $u^*$. Denote, now, the constrained optimal input with $u^*_c$, we have that $u^*_c \in \mathcal{H}_{++}$ and $H(u^*) < H(u^*_c)$. Being the Hamiltonian monotone with respect to the input, $H(u) \geq H(u^*_c) + \nabla_u H(u^*_c)^T(u - u^*_c) \geq H(u^*) + \nabla_u H(u^*)^T(u - u^*)$, i.e. there are not any values of $u$ which improves the Hamiltonian performance index. Hence \eqref{eq:sat_con} is an optimum for the input constrained problem.
\end{proof}

\subsection{Path Constraints}
In order to incorporate the path constraints within the OCP, we define an augmented cost function $\tilde{J}$ such that, when the state approaches the boundary of the forbidden set, its value becomes significantly larger than $J$, $\lim_{C(x,t)\to 0} \tilde{J} \gg J$. Therefore, we augment the cost function by choosing a suitable penalty function $k(x)$ which prevents the states approach the boundary of the constrained set whilst its value is (almost) null when far from the boundaries. This is a well known methodology to deal with the path constraints \cite{wang2014}.

\begin{assumption}
Assume the penalty function $k(x) \in \mathcal{C}^1 ( \mathbb{R}^{n_x} \rightarrow \mathbb{R})$.
\end{assumption}

The augmented cost index is given by
\begin{equation}\label{J_constr}
\begin{split}
\tilde{J}(u(\tau)) &= J(u(\tau)) + \int_{t}^{t + T_P}\sum_{i = 1}^n k_{i}(x) ~\mathrm{d}\tau .
\end{split}
\end{equation}
where $n$ is the number of the state constrains. The, the augmented Hamiltonian is
\begin{equation}
\begin{split}
\tilde{H}(x, u, \lambda) &= H(x,u,\lambda) + \sum_{i = 1}^n k_{i}(x).
\end{split}
\end{equation}
With the slight modification of the NMPC performance index and the consequent Hamiltonian augmentation, the contribute of the penalty function will affect the Euler-Lagrange equations by adding the terms of $\nabla_x \sum_{i = 1}^n k_{i}(x)$. In a more general form $\dot{\lambda} = -\nabla_x \bigl(H + \sum_{i = 1}^n k_{i}(x)\bigr)$. 

\section{Closed-loop Local Stability and Convergence}\label{finitegain}
The closed-loop stability for the nonlinear MPC schemes is an hard issue to tacke. According to \cite{chen1998}, the closed-loop stabilty for finite horizon can only be achieved by a suitable tuning of prediction horizon and weighting matrices. Over the last decades, important results, were obtained by \cite{mayne1990, mayne1991, mayne1993}, which posed the basis for future works on nonlinear systems stability. These latter works base the stability results on the differentiability and/or the Lipschitzianity of the optimal value function, as well as, an exact fulfillment of a terminal equality constraints, that, in the nonlinear case, is hardly satisfiable. This latter constraint has been relaxed in \cite{mayne1993} and \cite{chen1998}, being substituted by a terminal inequality constraint. Moreover, as pointed out by \cite{la2017} the main assumption proposed by the cited works is that the predicted trajectory coincides with the true trajectory. Conversely, we propose a stability criterion where this latter assumption does not hold.

In our work, the terminal equality constraints on the state is dropped, this is the common setup adopted in the recent years by \cite{la2017}, \cite{grimm2005}, and \cite{reble2012}. For the NMPC scheme at hand, the equality constraint has been substituted with the inequality constraint that bounds the state - at the end of the prediction horizon - in a prescribed terminal region. Furthermore, in \cite{mayne1990}, the final border of the prediction horizon is kept fixed at the time when the state reaches the equilibrium. Thanks to this assumption, when applying the receding horizon strategy, the trajectory between two different optimization steps will not change. However, this assumption does not always reflect the behavior of the real NMPC applications. Therefore, in our work, the final border of the prediction horizon will keep moving on, resulting in different state trajectory when considering different optimizations steps.


\begin{definition}\label{cl_loop}
The closed-loop system is defined as the system described by equation \eqref{eq:nonlinear}, where $u(t) = u^*(t_k)$, $t \in [t_k, t_{k+1}]$, $\forall t$, and $u^*$ is the solution of the optimization problem \eqref{eq:functional}, computed at each sampling time $t_k$.
\end{definition}

To this end, it is useful to introduce the concept of stability studied in this paper. This aspect is quite important since the proposed stability analysis is slightly different with respect to the classical ones in literature and it refers to a sort of practical stability.

\begin{definition}\label{as_stab_fin}
\textbf{Finite-time Practical Stability.}
Consider an autonomous nonlinear system, with state $x$. Let $x_r$ be an equilibrium point of the system and let $\epsilon > 0$. The set $\mathbf{B}(x_r, \epsilon)$ is locally finite-time stable if both the following conditions hold:\\
(i) It is locally stable. \\
(ii) It is locally attractive in finite-time: a $\delta_a > \epsilon$ and a finite $t_k \geq 0$ exist such that, for any initial condition $x(0) \in \mathbf{B}(x_r, \delta_a)$, it holds that $x(t) \rightarrow \mathbf{B}(x_r, \epsilon)$, $ \forall t \geq t_k$.
\end{definition}

Note that, the concept of stability introduced in the above definitions are similar (but slightly different) with respect to the classical ones available in the literature. Indeed, they refer to stability of a set containing the equilibrium point, rather than stability of the equilibrium point itself. According to Definition \ref{as_stab_fin}, practical stability requires that the system trajectory converges to set if the initial conditions are chosen sufficiently close to the set itself. It is similar to the classical concept of asymptotic stability of an equilibrium point by Lyapunov but, in the former case, the attractor is a set and not a single point.

\begin{assumption}\label{eq_isolated}
Assume that the nonlinear system $\dot{x} = f(x) + g(x)u$ has an isolated equilibrium point $(x_{eq}, u_{eq})$. 
\end{assumption}



\begin{assumption}\label{reachability}
Assume that a ball $\mathbf{B}(x_r, \epsilon_r)$ exists such that:
\begin{itemize}
\item The function $f(x) + g(x)u$ in \eqref{eq:nonlinear} is Lipschitz continuous with a Lipschitz constant $\Gamma$ for all $u \in U_C$ and all $x \in \mathbf{B}(x_r, \epsilon_r)$.
\item The constraints in the optimization problem \eqref{eq:functional} are feasible for all  $x \in \mathbf{B}(x_r, \epsilon_r)$
\item The reference is locally reachable for the given $T_p$. That is, for any $x(t) \in \mathbf{B}(x_r, \epsilon_r)$, a command signal $\hat{u}(t: t + T_p)$ exists such that $\hat{x}(t: t + T_p)$ satisfies the constraints in \eqref{eq:functional} and $\hat{x}(t + T_p) = x_r$.
\end{itemize}
\end{assumption}



\begin{theorem}\label{th:closedloop}\textbf{Closed-loop Local Finite-time Practical Stability}.
Consider the closed-loop system of Definition \ref{cl_loop} and the cost function defined in \eqref{eq:control_law_quad}. Let $x_r$ be an equilibrium point of this system. Let Assumptions \ref{lip}, 
 \ref{eq_isolated}, hold and Assumption \ref{reachability} for some $\epsilon_r > 0$. Then, for any $\epsilon \in (0,\epsilon_r]$ and any initial condition $x(0) \in \mathbf{B}(x_r, \delta_a)$, a finite diagonal matrix $\mathbf{P} > 0$ of \eqref{eq:control_law_quad} exists such that the ball $\mathbf{B}(x_r, \epsilon$) is finite-time practical stable.
\end{theorem}


\begin{proof}
{\it (Preliminary proof.)}
In order to make clearer to the reader the mathematical procedure for developing the closed-loop stability analysis, the proof is organized in different parts.\\
\\
\textbf{Part 1: Equivalent optimization problem}\\
Consider the optimization problem
\begin{equation}\label{eq:opt_equiv}
u^*(t_k:t_k + T_p) = \arg \min_{\hat{u}(\cdot),\eta_k} J_R\bigl(\hat{u}(t_k:t_k + T_p)\bigr)
\end{equation}
subject to all constraints in \eqref{eq:functional} and
\begin{equation}
\|\tilde{x}(t_k + T_p) \|_{P}^2 \leq \eta_k
\end{equation}
where
\begin{equation}
J_R\bigl(\hat{u}(t_k:t_k + T_p)\bigr) = \int_{t_k}^{t_k + T_p} \| \tilde{x}_p(\tau) \|^2_Q + \| \tilde{u}(\tau) \|^2_R ~\mathrm{d}\tau + \eta_k
\end{equation}
and $\tilde{u}(\tau) \doteq u_r - \hat{u}(\tau)$, $u_r \in \mathbb{R}^{n_u}$ is the reference input.  The terminal term of the cost function $\|\tilde{x}_p(t + T_p)\|_P^2$ in \ref{eq:control_law_quad} and \ref{eq:control_law_bang} can be re-written such that it can be included as an additional inequality constraint to the optimization problem. Namely, for a fixed matrix $\mathbf{P}>0$ and for each $k$, a $\eta_k \geq 0$ exists such that the problems \eqref{eq:functional} and \eqref{eq:opt_equiv} are equivalent. It follows that, for any $\epsilon_c \geq 0$, there exists a $\mathbf{P} > 0$ and a sequence $\eta_k \leq \epsilon_c$ for all $k$, such that the problems \eqref{eq:functional} and \eqref{eq:opt_equiv} are equivalent for all $k$. Furthermore, it holds that:
\begin{equation}\label{eq:constr_eq_mod}
\|\tilde{x}(t_k + T_p) \|_{P}^2 \leq \epsilon_c, ~~ \forall k
\end{equation}
In the following of the proof, the signal $\tilde{x}^*(t_k:t_k + T_p)$ obtained solving \eqref{eq:functional} (or \eqref{eq:opt_equiv} jointly with \eqref{eq:constr_eq_mod}) at time $k$ is denoted with $\tilde{x}_k^*(t_k:t_k + T_p)$ .\\
\\
\textbf{Part 2: Equivalent and backward systems}\\
The first equation of \eqref{eq:nonlinear} can rewritten in function of the relative state $\tilde{x} \doteq x_r - \hat{x}$, giving rise to the following equivalent system:
\begin{equation}\label{eq:tilde_sys}
\dot{\tilde{x}}  = \tilde{f}(\tilde{x}) + \tilde{g}(\tilde{x})u
\end{equation}
where $\tilde{f}(\tilde{x}) \doteq f(x_r - \hat{x})$, $\tilde{g}(\tilde{x}) \doteq g(x_r - \hat{x})$, and the disturbance has been supposed null. Define the backward time $t_b \doteq t_k + T_p - t$. For any given $k$, we have that $\dfrac{\mathrm{d}}{\mathrm{d}t_b}\tilde{x} = \dfrac{\mathrm{d}t}{\mathrm{d}t_b}\dfrac{\mathrm{d}}{\mathrm{d}t}\tilde{x} = -\dot{\tilde{x}}$. Thus, Equation \eqref{eq:back_sys} can be rewritten as
\begin{equation}\label{eq:back_sys}
\dfrac{\mathrm{d}}{\mathrm{d}t_b}\tilde{x} = - \bigl(\tilde{f}(\tilde{x}) + \tilde{g}(\tilde{x})u\bigr)
\end{equation}
which we call the backward system. \\
\\
\textbf{Part 3: Bound on} $\tilde{x}^*_{k}(t)$\\
The signal $\tilde{x}_k^*(t_k:t_k + T_p)$ is the solution to Equation \eqref{eq:tilde_sys}, corresponding to the initial condition $\tilde{x}_k^*(t_k) = \tilde{x}(t_k) = x_r - \hat{x}(t_k) = x_r - x(t_k)$ and the input signal $u^*(t_k:t_k + T_p)$. In a similar fashion, the signal $\tilde{x}^*_{bk}(0:T_p)$ is the solution to Equation \eqref{eq:back_sys}, corresponding to the initial condition $\tilde{x}^*_{bk}(0) = \tilde{x}(t_k + T_p)$ and the input signal $u^*(t_k + T_p:t_k)$. Clearly, it holds that $\tilde{x}^*_{bk}(t_b) = \tilde{x}^*_{k}(t_k + T_p - t_b)$, $\forall t_b \in [0, T_p]$. From \eqref{eq:constr_eq_mod}, we have that $\|\tilde{x}(t_k + T_p) \|_{P}^2 \leq \epsilon_c$. This implies $\|\tilde{x}^*_{bk}(0)\|_P^2 \leq \epsilon_c$, which in turn implies $\|\tilde{x}^*_{bk}(0)\|^2 \leq \dfrac{1}{p_m}\|\tilde{x}^*_{bk}(0)\|^2_P \leq \dfrac{\epsilon_c}{p_m}$, where $p_m \doteq \min_i p_i$ and $p_i, ~i= 1, ..., n_x$, are the diagonal elements of $\mathbf{P}$. From the Pontryagin theory (see Equation \eqref{eq:control_law_quad}), the command $u$ is continuous with respect to time. We can this apply the Theorem 2.1 (Chapter 1) stated in Coddington and Levinson \cite{coddington1955}, yielding following inequality:
\begin{equation}
\|\tilde{x}^*_{bk}(t_b)\| \leq \sqrt{\dfrac{\epsilon_c}{p_m}}e^{\Gamma t_b}
\end{equation}
where $\Gamma$ is the local Lipschitz constant of $\tilde{f} + \tilde{g}u$ at $0$ (which is equal to the local Lipschitz constant of $f + gu$ at $x_r$). Since $\tilde{x}^*_{bk}(t_b) = \tilde{x}^*_{k}(t_k + T_p - t_b)$, $\forall t_b \in [0, T_p]$ and $t_b \doteq t_k + T_p - t$, we have that the following inequalities chain:
\begin{equation}\label{eq:x_bound}
\|\tilde{x}^*_{bk}(t_b)\| \leq \sqrt{\dfrac{\epsilon_c}{p_m}}e^{\Gamma (t_k + T_p - t)} \leq \sqrt{\dfrac{\epsilon_c}{p_m}}e^{\Gamma T_p}, ~~ \forall k \geq 0.
\end{equation}
\\
\textbf{Part 4: Discrete-time Lyapunov-like function}\\
Define now the discrete Lyapunov-like function $V_k$ as
\begin{equation}\label{eq:lyap_def}
V_k \doteq \int_{t_k}^{t_k + T_p}\xi_k(\tau) ~ \mathrm{d}\tau.
\end{equation}
where $\xi_k(\tau) \doteq \|\tilde{x}^*_k(\tau) \|^2$, whereas $\tilde{x}^*_k(\tau)$ is the optimal state trajectory, obtained by solving \eqref{eq:opt_equiv} at time $t_k$. We have that
\begin{equation}\label{eq:lyap_var}
\begin{split}
V_{k+1} - V_k &= \int_{t_{k+1}}^{t_{k+1} + T_p}\xi_{k+1}(\tau) ~ \mathrm{d}\tau -  \int_{t_k}^{t_k + T_p}\xi_k(\tau) ~ \mathrm{d}\tau  =\\
&=  \int_{t_k + T_s}^{t_k + T_p + T_s}\xi_{k+1}(\tau) ~ \mathrm{d}\tau + \\
&- \int_{t_k + T_s}^{t_k + T_p}\xi_{k}(\tau) ~ \mathrm{d}\tau - \int_{t_k}^{t_k + T_s}\xi_{k+1}(\tau) ~ \mathrm{d}\tau =\\
&= - \int_{t_k}^{t_k + T_s}\xi_{k}(\tau) ~ \mathrm{d}\tau + \zeta_k
\end{split}
\end{equation}
where
\begin{equation}
\zeta_k \doteq \int_{t_k + T_s}^{t_k + T_p + T_s}\xi_{k+1}(\tau) ~ \mathrm{d}\tau - \int_{t_k + T_s}^{t_k + T_p}\xi_{k}(\tau) ~ \mathrm{d}\tau .
\end{equation}
\\
\textbf{Part 5: Bound on} $\zeta_k$\\
Since $\xi_k(\tau) \geq 0$, $\forall k$, $\forall \tau$, the following inequality hold:
\begin{equation}
 - \int_{t_k + T_s}^{t_k + T_p}\xi_{k}(\tau) \leq \zeta_k \leq \int_{t_k + T_s}^{t_k + T_p + T_s}\xi_{k+1}(\tau).
\end{equation}
Using \eqref{eq:x_bound} for the time intervals $[t_k + T_S, t_k + T_p]$ and $[t_k + T_S, t_k + T_p + T_S]$, we obtain that
\begin{equation}
\int_{t_k + T_s}^{t_k + T_p + T_s}\xi_{k+1}(\tau) ~ \mathrm{d}\tau \leq \dfrac{\epsilon_cT_p}{p_m}e^{2\Gamma T_p}
\end{equation}
\begin{equation}
\int_{t_k + T_s}^{t_k + T_p}\xi_{k}(\tau) ~ \mathrm{d}\tau \leq \dfrac{\epsilon_c(T_p-T_S)}{p_m}e^{2\Gamma (T_p- T_S)}.
\end{equation}
It follows that $\zeta_k$ is bounded as
\begin{equation}
|\zeta_k| \leq \dfrac{\epsilon_cT_p}{p_m}e^{2\Gamma T_p} \doteq \epsilon_c \bar{\nu}
\end{equation}
\\
\textbf{Part 6: Time evolution of} $V_k$\\
Note that $\int_{t_k}^{t_k + T_S}\xi_{k}(\tau) ~ \mathrm{d}\tau$ is a fraction of $V_k$. That is,  $\int_{t_k}^{t_k + T_S}\xi_{k}(\tau) ~ \mathrm{d}\tau = \phi_k V_k$ for some $0 \leq \phi_k \leq 1$. Equation \eqref{eq:lyap_var} can this be re-written as
\begin{equation}
V_{k+1} = a_kV_k + \zeta_k
\end{equation}
where $0 \leq a_k \doteq 1 - \phi_k \leq 1, ~\forall k$. This equation describes a linear time-varying system with input $\zeta_k$ and state $V_k$. From linear system theory (see, e.g., Rugh \cite{rugh1996}), the solution of this system is given by
\begin{equation}\label{eq:lyap_sys}
V_k = \Phi_{k,0}V_0 + \sum_{l = 0}^{k-1}\Phi_{k,l+1}\zeta_l
\end{equation}
where $\Phi_{k,l} \doteq \prod_{i = l}^{k-1}a_i, ~ k>l \geq 0, ~\Phi_{k,k} \doteq 1$. We can now distinguish two alternative cases.\\
\underline{Case 1:}
A $k\geq 0 $ exists such that $a_k = 1$ or, equivalently, $\phi_k = 0$. Then, $\int_{t_k}^{t_k + T_S}\xi_{k}(\tau) = 0$, which in turn implies that $\xi_{k}(\tau) \doteq \|\tilde{x}_k^*(\tau)\|^2 = 0, ~ \forall \tau \in [t_k, t_k + T_S]$. Observing that $\tilde{x}_k^*(t_k) \doteq x_r - x(t_k)$, the equality $ \|\tilde{x}_k^*(\tau)\|^2 = 0$ means that $x(t_k) = x_r$. Since $\mathbf{Q}$, $\mathbf{R}$, and $\mathbf{P}$ are positive definite, solving the optimization problem \eqref{eq:opt_equiv} (or \eqref{eq:functional}) with this initial condition gives $u^*(t_k:t_k+T_p) = u_r$. Indeed, with $x(t_k) = x_r$ and $u^*(t_k:t_k+T_p) = u_r$, we have that $J \bigl(u^*(t_k:t_k+T_p)\bigr) = J_R \bigl(u^*(t_k:t_k+T_p)\bigr) = 0$. These considerations show that $x_r$ is an equilibrium of the closed-loop system, implying that $x(t_j) = x_r, ~\forall j \geq k$.\\
\underline{Case 2:}
A $\bar{a} < 1$ exists such that $a_k \leq \bar{a}, ~\forall k$. Consider that
\begin{equation}\label{eq:lyap_bound}
\begin{split}
 \sum_{l = 0}^{k-1}\Phi_{k,l+1}\zeta_l &\leq \epsilon_c \bar{\nu} \sum_{l = 0}^{k-1}\bar{a}^{k-l-1} = \epsilon_c \bar{\nu}\bar{a}^{k-1} \sum_{l = 0}^{k-1}\bar{a}^{-l}  \\
 &=\epsilon_c \bar{\nu}\bar{a}^{k-1} \sum_{l = 0}^{k-1}\biggl(\dfrac{1}{\bar{a}}\biggr)^l = \\
 &=\epsilon_c \bar{\nu}\bar{a}^{k-1}\dfrac{1 - 1/\bar{a}^k}{1 - 1/\bar{a}} = \epsilon_c \bar{\nu}\dfrac{\bar{a}^{k-1} - 1/\bar{a}}{1 - 1/\bar{a}} = \\
 &=\epsilon_c \bar{\nu}\dfrac{\bar{a}^{k} - 1}{\bar{a}-1} = \epsilon_c \bar{\nu}\dfrac{1 -\bar{a}^{k}}{1 - \bar{a}} \leq \dfrac{\epsilon_c \bar{\nu}}{1 - \bar{a}}, ~~\forall k
\end{split}
\end{equation}
where the geometric series formula has been used. Morevoer, $\Phi_{k,0} \rightarrow 0$ as $k \rightarrow \infty$. From \eqref{eq:lyap_sys} and \eqref{eq:lyap_bound}, we thus obtain the following asymptotic inequality:
\begin{equation}
\lim_{k \rightarrow \infty}V_k \leq \dfrac{\epsilon_c \bar{\nu}}{1 - \bar{a}} .
\end{equation}
For every $k$ and $\tilde{x}_k^*(t_k) \neq x_r$, $V_k$ is a locally positive-definite function of $\tilde{x}_k^*(t_k)$. It follows that, for any $\epsilon_x > 0$, an $\epsilon_c$ and a finite $k \geq 0$ exists, such that $\| \tilde{x}_j(t_j)\| = \| x_r - x(t_j)\| \leq \epsilon_x$, $\forall j \geq k$.\\
\\
\\
\\
\textbf{Part 7: Proof of Practical Finite-time Stability}\\
From the study of the two cases in the Part 6 of the proof, jointly with the consideration of the Part 1, we have that, for any $\epsilon_x > 0$, a matrix $\mathbf{P} >0$ and a finite $k \geq 0$ exists such that $\|\tilde{x}_j^*(t_j) \| = \| x_r - x(t_j)\| \leq \epsilon_x, ~\forall j \geq k$. Hence, from Theorem 2.1 (Chapter 1) stated in Coddington and Levinson \cite{coddington1955}, we obtain the following bound:
\begin{equation}
\| x_r - x(t) \| \leq \epsilon_x e^{\Gamma(t-t_j)} \leq \epsilon_x e^{\Gamma T_S} \doteq \epsilon, ~~\forall t \in [t_j, t_{j+1}], ~~\forall j \geq k .
\end{equation}
We can conclude that, for any $\epsilon > 0$, a matrix $\mathbf{P} > 0$ and a finite $k \geq 0$ exist such that $\| x_r - x(t) \| \leq \epsilon$, $\forall t \geq t_k$, which is the set finite-time practical stability definition.


\end{proof}

\section{SIMULATED EXAMPLES}\label{sim}
Consider the predatory-prey Lotka-Volterra model, described by a couple of first-order nonlinear differential equations with an exogeneous input applied on both states:
\begin{equation}\label{eq:vanderpol}
\begin{cases}
\dot{x}_1 = x_1(\alpha - \beta x_2) + x_1u_1  \\
\dot{x}_2 = x_2(\gamma x_1 - \delta) + x_2u_2 \\
\end{cases}
\end{equation}
where $x_1$ and $x_2$ are the prey and predator population respectively and $u_1$ and $u_2$ the corresponding input components. Let $\alpha = 0.25$, $\beta = 0.25$, $\gamma = 0.008$, and $\delta = 0.008$ be parameters describing the interaction between the two species. The admissible input set is described by $U_C = \{u(t) \in \mathbb{R}: ~ -u_{i_{max}} \leq u_i(t)\leq u_{i_{max}}, ~\forall t \}$, where $u_{1_{max}}=10$ and $u_{2_{max}}=5$. Concerning the state constraints, we designed a nonlinear function which prevents the predator specie grows too abruptly with respect to the prey specie, then, avoiding the extinction of both species when the prey population goes to zero.
Hence, $X_C = \{x(t) \in \mathbb{R}^2: 5 - \bigl((x_1 - 100)^2 + (x_2 - 51.5)^2\bigr)^{1/2} \leq 0, ~\forall t\}$. Thus, the state constraints are handled employing a Gaussian-like penalty function $k(x) = a \exp{(-bC(x)^2)}$ with $C =  5 - \bigl((x_1 - 100)^2 + (x_2 - 51.5)^2\bigr)^{1/2}$, $a = 10^6$, and $b = 1$. The augmented Hamiltonian is
\begin{equation}
\begin{split}
\tilde{H} &=  \lambda_1\bigl( x_1(\alpha - \beta x_2) + x_1u_1\bigr)+ \lambda_2\bigl(x_2(\gamma x_1 - \delta) + x_2u_2\bigr)
+ \\ 
&+ \sum_{i}^{n_u}\mathbf{R}_{i,i}u_i^2 + \sum_{i}^{n_x}\mathbf{Q}_{i,i}\tilde{x}_{p_i}^2 + k(x) 
\end{split}
\end{equation}
and the TPBVP is formalized as:
\begin{equation}\label{eq:BVP_PC}
\begin{cases}
\dot{x}_1 &= x_1(\alpha - \beta x_2) + x_1u_1  \\
\dot{x}_2 &= x_2(\gamma x_1 - \delta) + x_2u_2  \\
\dot{\lambda}_1 &= -\alpha\lambda_2 + \beta\lambda_1\lambda_2 - \lambda_1u_1 -\gamma\lambda_2x_2  - 2\mathbf{Q}_{11}\tilde{x}_{p_1} - \dfrac{\partial{k(x)}}{\partial{x_1}} \\
\dot{\lambda}_2 &=\beta\lambda_1x_1 - \gamma\lambda_2x_1 + \delta\lambda_2 - \lambda_2u_2 - 2\mathbf{Q}_{22}\tilde{x}_{p_2} - \dfrac{\partial{k(x)}}{\partial{x_2}}\\
x_i &= x(t_k)  \\
\lambda(t_f) &= \bigl(2\mathbf{P}\tilde{x}_p(t_f)\bigr)^T
\end{cases}
\end{equation}
The solution of the TPBVP in \eqref{eq:BVP_PC} provides the $\lambda$ for the explicit optimal control law:
\begin{equation}
u^* = -\dfrac{1}{2}\mathbf{R}^{-1}\lambda^Tx
\end{equation}
The NMPC parameters are listed in Table \ref{tb:nmpc_LV}.
\begin{table}
\caption{NMPC Parameters}\label{tb:nmpc_LV}
\begin{center}
\begin{tabular}{|c|c|c|c|c|}
\hline
{$T_S$}&{$T_p$}&{$\textbf{R}$}&{$\textbf{Q}$}&{$\textbf{P}$}\\
\hline
{$0.0001~s$}&{$0.001~s$}&{$500\cdot\mathbf{I}_{2\times 2}$}&{$diag(10, 35)$}&{$diag(10, 35)$} \\
\hline
\end{tabular}
\end{center}
\end{table}
The desired state is $[x_{r_1}(t),x_{r_2}(t)]^T=[10\cos(t)+100,10\sin(t)+50]^T$. The initial state is $x_0=[40,40]^T$. This means that prey and predator populations are far more than zero, that is the two species are both far from the risk of extinction.

\begin{figure}[thpb]
      \centering
      \includegraphics[scale=0.6]{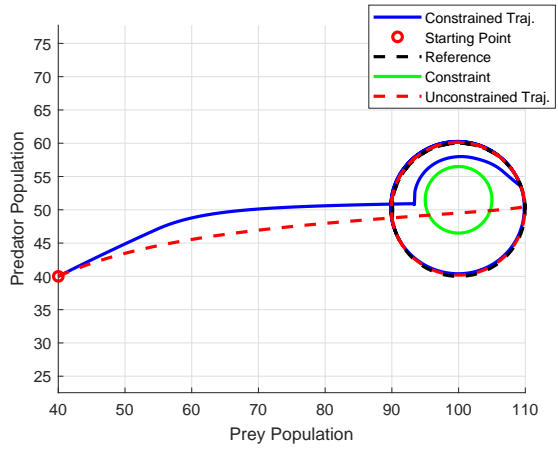}
      \caption{Predator-prey populations phase-plane}
\label{fig:LV_trajec}
   \end{figure}

\begin{figure}[thpb]
      \centering
      \includegraphics[scale=0.5]{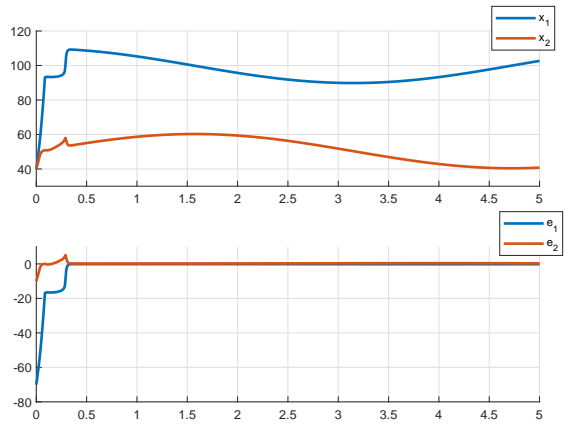}
   \caption{Temporal evolution of populations and corresponding tracking errors}
\label{fig:LV_state}
   \end{figure}

\begin{figure}[thpb]
      \centering
      \includegraphics[scale=0.6]{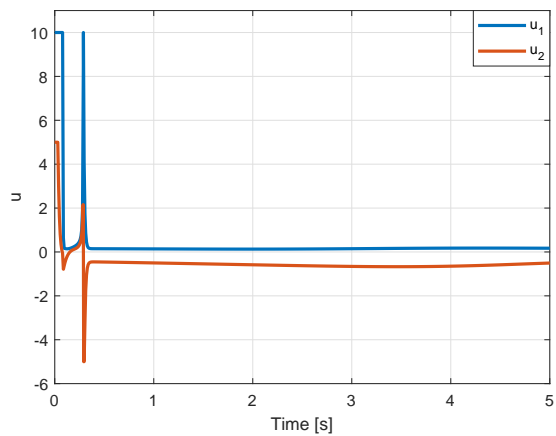}    
\caption{Control components}
\label{fig:LV_command}
   \end{figure}

In Figure \ref{fig:LV_trajec}, the phase-plane curve of predator-prey populations is shown. In particular, it is highlighted how the NMPC approach is perfecty able to fulfill the input and state constraints, without affecting the tracking performance. Figure \ref{fig:LV_state} displays the time evolution of populations $x_1$ and $x_2$, and the corresponding tracking errors $e_1$ and $e_2$. This latter have a very fast convergence to zero, proving the effectiveness of the optimization algorithm. Finally, in Figure \ref{fig:LV_command}, the command activity is reported.

\begin{figure}[thpb]
      \centering
      \includegraphics[scale=0.55]{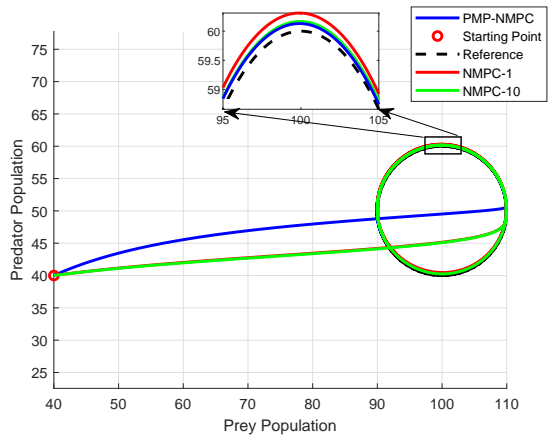}    
\caption{PMP-NMPC and Numerical NMPCs phase-plane}
\label{fig:NMPC_comparison}
   \end{figure}

We are now interested in comparing the behavior of the solutions when employing different optimization strategies: the PMP-based and the Sequential Quadratic Programming (SQP) solutions. Concerning the SQP case, we further considered two different cases: i) constant input parametrization (NMPC-1), ii) piece-wise constant input parametrization with $N=10$ (NMPC-10). In the latter case, the input is parametrized with the same sampling steps adopted in the PMP-based solution.
Figure \ref{fig:NMPC_comparison} reports the results obtained in the unconstrained case, both for SQP and PMP. The resulting trajectories are slightly different. However, in all configurations, the NMPC is able to get a good tracking of the reference.

\begin{figure}[thpb]
      \centering
      \includegraphics[scale=0.6]{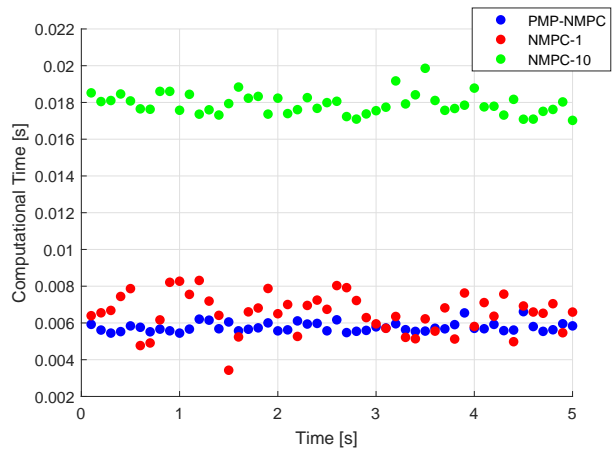}    
\caption{Computational Cost analysis}
\label{fig:LV_comp_time}
   \end{figure}

From the computational burden point of view, Figure \ref{fig:LV_comp_time} presents a comparison between the solutions. If considering a similar input parametrization, the PMP-NMPC shows superior computational performances with respect to the SQP-NMPC-10. Moreover, also when considering the constant input parametrization, the PMP-NMPC owns slight better performances - together with a better reference tracking - with respect to the SQP solution.

Then, we highlight that the advantages of the proposed PMP-based NMPC framework are: i) a better reference tracking than the NMPC-10 configuration, ii) a similar computational cost with respect to the NMPC-1 configuration.
\section{CONCLUSIONS}\label{conclusion}
We proposed an alternative approach for the Nonlinear Model Predictive Control optimization problem. We obtained a control law by developing an algorithm based on the Pontryagin Minimum Principle turning the optimal control problem into a two-points boundary value problem. The resulting optimal input is function of the co-state variables, whose time evolution is described by the Euler-Lagrange differential equation. Hence, the optimal control law was obtained analytically by minimizing the Hamiltonian of the system. Moreover, we also coped with state constraints by exploiting a suitable penalty function within the cost function without any modification of the optimization algorithm. The proposed methodology was then applied to the Lokta-Volterra dynamics. The results highlighted the effectiveness of the control algorithm, showing excellent reference tracking and the compliance with the input and path constraints.

\end{document}